\documentclass[11pt,reqno]{amsart}
\usepackage{graphicx}
\usepackage{verbatim}
\usepackage{textcomp}
\usepackage{amssymb}
\usepackage{cite}
\usepackage{amsmath}
\usepackage{latexsym}
\usepackage{amscd}
\usepackage{amsthm}
\usepackage{mathrsfs}
\usepackage{xypic}
\usepackage{bm}
\usepackage{url}
\usepackage{hyperref}

\newtheorem{thm}{Theorem}[section]
\newtheorem{corr}[thm]{Corollary}

\theoremstyle{definition}

\theoremstyle{remark}
\newtheorem{rem}{Remark}[section]
\numberwithin{equation}{section}
\begin{document}
\title[Gradient estimates and Liouville type theorems]
{Gradient estimates and Liouville type theorems for a nonlinear
elliptic equation}
\author{Guangyue Huang}
\address{Department of Mathematics, Henan Normal
University, Xinxiang 453007, P.R. China} \email{hgy@henannu.edu.cn }
\author{Bingqing Ma}
\address{Department of Mathematics, Henan Normal
University, Xinxiang 453007, P.R. China} \email{bqma@henannu.edu.cn}
\thanks{The research of the first
author is supported by NSFC(No. 11371018, 11171091). The research of
the second author is supported by NSFC (No. 11401179).} \maketitle

\begin{abstract} Let $(M^n,g)$ be an $n$-dimensional
complete Riemannian manifold. We consider gradient estimates and
Liouville type theorems for positive solutions to the following
nonlinear elliptic equation:
$$\Delta u+au\log u=0,$$ where $a$ is a nonzero constant. In particular, for $a<0$,
we prove that any bounded positive solution of the above equation
with a suitable condition for $a$ with respect to the lower bound of
Ricci curvature must be $u\equiv 1$. This generalizes a classical
result of Yau.

\end{abstract}

\section{Introduction}

In this paper, we study positive solutions of the equation
\begin{equation}\label{Int1}
\Delta u+au\log u=0\end{equation} on an $n$-dimensional complete
Riemannian manifold $(M^n,g)$, where $a$ is a nonzero constant. The
first study of this and related nonlinear equations can be traced
back to Li \cite{Li1991}, and later by Ma \cite{Ma06} and Yang
\cite{Yangyuyan07}, who derived various gradient estimates and
Harnack estimates and noted the relation to gradient Ricci solitons,
which are self-similar solutions to Ricci flow and arise in the
blow-up procedure of the long time existence or convergence of the
flow. Moreover, the equation \eqref{Int1} is closely related to the
famous Gross Logarithmic Sobolev inequality, see \cite{Gross93}. For
the recent research of  \eqref{Int1}, on can refer to
\cite{Cao2013,CC2009,HHL2013,Wu2010,Qian10} and the references
therein.

It is well-known that for complete noncompact Riemannian manifolds
with ${\rm Ric}\geq0$, Yau \cite{yau75} showed that every positive
or bounded solution to the equation
\begin{equation}\label{Int2} \Delta u=0\end{equation}
is constant. We note that the equation \eqref{Int2} can be seen as a
special case of \eqref{Int1} when $a=0$. Therefore, a natural idea
is to achieve similar Liouville type theorems for positive solutions
to the nonlinear elliptic equation \eqref{Int1}.

Inspired by the method used by Brighton in \cite{Brighton2013}, this
paper is concerned with gradient estimates and Liouville type
theorems for positive solutions to the nonlinear elliptic equation
\eqref{Int1} with $a\neq0$ and obtain the following results:

\begin{thm}\label{thm1-1}
Let $(M^n,g)$ be an $n$-dimensional complete Riemannian manifold
with ${\rm Ric}(B_p(2R))  \geq-K$, where $K\geq0$ is a constant.
Suppose that $u$ is a positive solution to \eqref{Int1} on
$B_p(2R)$. Then on $B_p(R)$, the following inequalities hold:

(1) If $a>0$, then
\begin{equation}\label{1thm1}\aligned  |\nabla u|
\leq&M\sqrt{\frac{(n+3)^2}{2n}(a+K)+\frac{C}{R^2}\left( 1+ \sqrt{K}R
\coth (\sqrt{K} R)\right)};
\endaligned\end{equation}

(2) If $a<0$, then
\begin{equation}\label{1thm2}\aligned  |\nabla u|
\leq&M\sqrt{\frac{(5n+6)^2}{36n}\max\{0,\ a+K \}+\frac{C}{R^2}\left(
1+ \sqrt{K}R \coth (\sqrt{K} R)\right)},
\endaligned\end{equation}
where $M=\sup\limits_{x\in B_p(2R)}u(x)$ and the constant $C$
depends only on $n$.

\end{thm}

Letting $R\rightarrow \infty$, we obtain the following gradient
estimates on complete noncompact Riemannian manifolds:

\begin{corr}\label{corr1-1}
Let $(M^n,g)$ be an $n$-dimensional complete noncompact Riemannian
manifold with ${\rm Ric}\geq-K$, where $K\geq0$ is a constant.
Suppose that $u$ is a positive solution to \eqref{Int1}. Then the
following inequalities hold:

(1) If $a>0$, then
\begin{equation}\label{1corr1}\aligned  |\nabla u|
\leq&\frac{(n+3)M}{\sqrt{2n}}\sqrt{a+K};
\endaligned\end{equation}

(2) If $a<0$, then
\begin{equation}\label{1corr2}\aligned  |\nabla u|
\leq&\frac{(5n+6)M}{6\sqrt{n}}\sqrt{\max\{0,\,a+K \}},
\endaligned\end{equation}
where $M=\sup\limits_{x\in M^n}u(x)$.

\end{corr}

In particular, for $a<0$, if $a\leq-K$, then \eqref{1corr2} shows
that $\max\{0,\, a+K \}=0$ and  $|\nabla u| \leq0$ whenever $u$ is a
bounded positive solution to \eqref{Int1}. This implies that
$u\equiv 1$ is a constant. Therefore, the following Liouville-type
result follows:

\begin{corr}\label{corr1-2}
Let $(M^n,g)$ be an $n$-dimensional complete noncompact Riemannian
manifold with ${\rm Ric}\geq-K$, where $K\geq0$ is a constant.
Suppose that $u$ is a bounded positive solution defined on $M^n$ to
\eqref{Int1} with $a<0$. If $a\leq-K$, then $u\equiv 1$ is a
constant.

\end{corr}

In particular, we have the following result:

\begin{corr}\label{corr1-3}
Let $(M^n,g)$ be an $n$-dimensional complete noncompact Riemannian
manifold with ${\rm Ric}\geq0$. Suppose that $u$ is a bounded
positive solution defined on $M^n$ to \eqref{Int1} with $a<0$, then
$u\equiv 1$ is a constant.

\end{corr}

\begin{rem}\label{rem1}
Clearly, our Corollary \ref{corr1-3} generalizes a uniqueness result
of Yau with respect to the heat equation to the nonlinear elliptic
equation \eqref{Int1} with $a<0$.
\end{rem}

\begin{rem}\label{rem2}
For an $n$-dimensional complete noncompact Riemannian manifold with
${\rm Ric}\geq0$, it has been shown in Corollary 1.1 of
\cite{HHL2013} that for positive solutions of \eqref{Int1} with
$a<0$,
\begin{equation}\label{rem1}
\frac{|\nabla u|^2}{u^2}+a\alpha\log u \leq \frac{-na\alpha^2}{8}
\end{equation} for any $\alpha>1$. From \eqref{rem1}, we can not obtain the uniqueness
theorem for any bounded positive solution of \eqref{Int1}.
Therefore, our Corollary \ref{corr1-3} generalizes (2) in Corollary
1.1 of \cite{HHL2013} in this sense.

\end{rem}

\section{Proof of results}

The proof of main results will follow from applying the Bochner
formula to an appropriate function $h$ of a given positive solution
$u$.

Let $h=u^\epsilon$, where $\epsilon\neq0$ is a constant to be
determined.  Then we have
\begin{equation}\label{Proof1}
\log h=\epsilon\log u.\end{equation} A straightforward computation
gives
\begin{equation}\label{Proof2}\aligned
\Delta
h=&\Delta(u^\epsilon)=\epsilon(\epsilon-1)u^{\epsilon-2}|\nabla
u|^2+\epsilon u^{\epsilon-1}\Delta u\\
=&\epsilon(\epsilon-1)u^{\epsilon-2}|\nabla u|^2-a\epsilon u^{\epsilon}\log u\\
=&\frac{\epsilon-1}{\epsilon} \frac{|\nabla h|^2}{h}-ah\log h.
\endaligned\end{equation}
Hence, we have
\begin{equation}\label{Proof3}\aligned
\nabla h\nabla\Delta h=&\nabla
h\nabla\Big(\frac{\epsilon-1}{\epsilon}
\frac{|\nabla h|^2}{h}-ah\log h\Big)\\
=&\frac{\epsilon-1}{\epsilon}\nabla h\nabla \frac{|\nabla
h|^2}{h}-a\nabla h\nabla(h\log h)\\
=&\frac{\epsilon-1}{\epsilon h}\nabla h\nabla (|\nabla
h|^2)-\frac{\epsilon-1}{\epsilon}\frac{|\nabla h|^4}{h^2}-ah\log
h\frac{|\nabla h|^2}{h}-a|\nabla h|^2.
\endaligned\end{equation}
Applying \eqref{Proof2} and \eqref{Proof3} into the well-known
Bochner formula to $h$, we have
\begin{equation}\label{Proof4}\aligned
\frac{1}{2}\Delta|\nabla h|^2=&|\nabla^2h|^2+\nabla h\nabla\Delta
h+{\rm Ric} (\nabla h,\nabla
h)\\
\geq&\frac{1}{n}(\Delta h)^2+\nabla h\nabla\Delta
h-K|\nabla h|^2\\
=&\frac{1}{n}\Big(\frac{\epsilon-1}{\epsilon} \frac{|\nabla
h|^2}{h}-ah\log h\Big)^2+\frac{\epsilon-1}{\epsilon }\frac{\nabla
h}{h}\nabla
(|\nabla h|^2)\\
&-\frac{\epsilon-1}{\epsilon}\frac{|\nabla h|^4}{h^2}
-ah\log h\frac{|\nabla h|^2}{h}-(a+K)|\nabla h|^2\\
=&\Big(\frac{(\epsilon-1)^2}{n\epsilon^2}-\frac{\epsilon-1}{\epsilon}\Big)\frac{|\nabla
h|^4}{h^2}-a\Big(\frac{2(\epsilon-1)}{n\epsilon}+1\Big)h\log
h\frac{|\nabla
h|^2}{h}\\
&+\frac{a^2}{n}(h\log h)^2+\frac{\epsilon-1}{\epsilon }\frac{\nabla
h}{h}\nabla (|\nabla h|^2)-(a+K)|\nabla h|^2.
\endaligned\end{equation}

Now we let
\begin{equation}\label{Proof5}
a\Big(\frac{2(\epsilon-1)}{n\epsilon}+1\Big)\geq0.
\end{equation}
Then for a fixed point $p$, if there exist a positive constant
$\delta$ such that $h\log h\leq\delta\frac{|\nabla h|^2}{h}$,
\eqref{Proof4} becomes
\begin{equation}\label{Proof6}\aligned
\frac{1}{2}\Delta|\nabla
h|^2\geq&\Big[\frac{(\epsilon-1)^2}{n\epsilon^2}-\frac{\epsilon-1}{\epsilon}
-a\delta\Big(\frac{2(\epsilon-1)}{n\epsilon}+1\Big)\Big]\frac{|\nabla
h|^4}{h^2}\\
&+\frac{a^2}{n}(h\log h)^2+\frac{\epsilon-1}{\epsilon }\frac{\nabla
h}{h}\nabla (|\nabla h|^2)-(a+K)|\nabla h|^2\\
\geq&\Big[\frac{(\epsilon-1)^2}{n\epsilon^2}-\frac{\epsilon-1}{\epsilon}
-a\delta\Big(\frac{2(\epsilon-1)}{n\epsilon}+1\Big)\Big]\frac{|\nabla
h|^4}{h^2}\\
&+\frac{\epsilon-1}{\epsilon }\frac{\nabla h}{h}\nabla (|\nabla
h|^2)-(a+K)|\nabla h|^2.
\endaligned\end{equation}
On the contrary, at the point $p$, if $h\log
h\geq\delta\frac{|\nabla h|^2}{h}$, then \eqref{Proof4} becomes
\begin{equation}\label{Proof7}\aligned
\frac{1}{2}\Delta|\nabla
h|^2\geq&\Big(\frac{(\epsilon-1)^2}{n\epsilon^2}-\frac{\epsilon-1}{\epsilon}\Big)\frac{|\nabla
h|^4}{h^2}
+\Big[\frac{a^2}{n}-\frac{a}{\delta}\Big(\frac{2(\epsilon-1)}{n\epsilon}+1\Big)\Big](h\log
h)^2\\
&+\frac{\epsilon-1}{\epsilon }\frac{\nabla h}{h}\nabla (|\nabla
h|^2)-(a+K)|\nabla h|^2\\
\geq&\Big\{\Big(\frac{(\epsilon-1)^2}{n\epsilon^2}-\frac{\epsilon-1}{\epsilon}\Big)+\delta^2\Big[\frac{a^2}{n}
-\frac{a}{\delta}\Big(\frac{2(\epsilon-1)}{n\epsilon}+1\Big)\Big]\Big\}\frac{|\nabla
h|^4}{h^2}\\
&+\frac{\epsilon-1}{\epsilon }\frac{\nabla h}{h}\nabla (|\nabla
h|^2)-(a+K)|\nabla h|^2\\
\geq&\Big\{\Big(\frac{(\epsilon-1)^2}{n\epsilon^2}-\frac{\epsilon-1}{\epsilon}\Big)
-a\delta\Big(\frac{2(\epsilon-1)}{n\epsilon}+1\Big)\Big\}\frac{|\nabla
h|^4}{h^2}\\
&+\frac{\epsilon-1}{\epsilon }\frac{\nabla h}{h}\nabla (|\nabla
h|^2)-(a+K)|\nabla h|^2
\endaligned\end{equation} as long as
\begin{equation}\label{Proof8}
\frac{a^2}{n}-\frac{a}{\delta}\Big(\frac{2(\epsilon-1)}{n\epsilon}+1\Big)>0.
\end{equation}

In order to obtain the bound of $|\nabla h|$ by using the maximum
principle for \eqref{Proof7}, it is sufficient to choose the
coefficient of $\frac{|\nabla h|^4}{h^2}$ in \eqref{Proof7} is
positive, that is,
\begin{equation}\label{Proof9}
\Big(\frac{(\epsilon-1)^2}{n\epsilon^2}-\frac{\epsilon-1}{\epsilon}\Big)
-a\delta\Big(\frac{2(\epsilon-1)}{n\epsilon}+1\Big)>0.\end{equation}

We next consider two cases:

{\bf Case one:} $a>0$.

In this case, provided
$\epsilon\in(\frac{2}{n+2},\frac{6}{(5-\sqrt{13})n+6})$, there will
exist an $\delta$ satisfying \eqref{Proof5}, \eqref{Proof8} and
\eqref{Proof9}. In particular, we choose
\begin{equation}\label{Proof10}
\epsilon=\frac{3}{n+3}
\end{equation} and
\begin{equation}\label{Proof11}
\delta=\frac{n}{2a},
\end{equation} then
\eqref{Proof7} becomes
\begin{equation}\label{Proof12}\aligned
\frac{1}{2}\Delta|\nabla h|^2\geq&\frac{5n}{18}\frac{|\nabla
h|^4}{h^2}-\frac{n}{3 }\frac{\nabla h}{h}\nabla (|\nabla
h|^2)-(a+K)|\nabla h|^2.
\endaligned\end{equation}

{\bf Case two:} $a<0$.

In this case, provided
$\epsilon\in(\frac{6}{(5+\sqrt{13})n+6},\frac{2}{n+2})$, there will
exist an $\delta$ satisfying \eqref{Proof5}, \eqref{Proof8} and
\eqref{Proof9}. In particular, we choose
\begin{equation}\label{Proof13}
\epsilon=\frac{6}{5n+6}
\end{equation} and
\begin{equation}\label{Proof14}
\delta=-\frac{3n}{4a}.
\end{equation} In particular, in this
case, \eqref{Proof7} becomes
\begin{equation}\label{Proof15}\aligned
\frac{1}{2}\Delta|\nabla h|^2\geq&\frac{37n}{36}\frac{|\nabla
h|^4}{h^2}-\frac{5n}{6}\frac{\nabla h}{h}\nabla (|\nabla
h|^2)-(a+K)|\nabla h|^2.
\endaligned\end{equation}

\vspace*{2mm}

 \noindent{\bf Proof of Theorem \ref{thm1-1}.}

First, we prove the case of $a>0$.

Denote by $B_p(R)$ the geodesic ball centered at $p$ with radius
$R$. Take a cut-off function $\phi$ (see \cite{schoenyau})
satisfying $\mathrm{supp}(\phi) \subset B_p(2R)$, $\phi |_{B_p(R)}
=1$ and
\begin{equation}\label{Proof16}
\frac{|\nabla \phi|^2}{\phi}\leq \frac{C}{R^2},\end{equation}
\begin{equation}\label{Proof17}-\Delta \phi \leq\frac{C}{R^2}\left( 1+ \sqrt{K}R \coth (\sqrt{K}
R)\right),\end{equation} where $C$ is a constant depending only on
$n$. Define $G=\phi |\nabla h|^2$. Next we will apply maximum
principle to $G$ on $B_p(2 R)$. Assume $G$ achieves its maximum at
the point $x_0\in B_p(2R)$ and assume $G(x_0)>0$ (otherwise the
proof is trivial). Then at the point $x_0$, it holds that
$$\Delta G\leq 0, \ \ \
\nabla (|\nabla h|^2)=-\frac{|\nabla h|^2}{\phi} \nabla \phi$$ and
\begin{equation}\label{Proof18}\aligned
0\geq& \Delta G\\
=&\phi\Delta(|\nabla h|^2)+|\nabla h|^2\Delta\phi+2\nabla\phi\nabla
|\nabla h|^2 \\
=&\phi\Delta(|\nabla h|^2)+\frac{\Delta\phi}{\phi}G
-2\frac{|\nabla\phi|^2}{\phi^2}G\\
\geq&2\phi\Big[\frac{5n}{18}\frac{|\nabla h|^4}{h^2}-\frac{n}{3
}\frac{\nabla h}{h}\nabla (|\nabla h|^2)-(a+K)|\nabla
h|^2\Big]\\
&+\frac{\Delta\phi}{\phi}G -2\frac{|\nabla\phi|^2}{\phi^2}G\\
=&\frac{5n}{9}\frac{G^2}{\phi
h^2}+\frac{2nG}{3\phi}\nabla\phi\frac{\nabla
h}{h}-2(a+K)G\\
&+\frac{\Delta\phi}{\phi}G -2\frac{|\nabla\phi|^2}{\phi^2}G,
\endaligned\end{equation} where the second inequality used
\eqref{Proof12}. Multiplying both sides of \eqref{Proof18} by
$\frac{\phi}{G}$ yields
\begin{equation}\label{Proof19}\aligned
\frac{5n}{9}\frac{G}{ h^2}\leq-\frac{2n}{3}\nabla\phi\frac{\nabla
h}{h}+2\phi(a+K)- \Delta\phi+2\frac{|\nabla\phi|^2}{\phi}.
\endaligned\end{equation}
Applying the Cauchy inequality
$$\aligned-\frac{2n}{3}\nabla\phi\frac{\nabla
h}{h}\leq&\frac{2n}{3}|\nabla\phi|\frac{|\nabla
h|}{h}\\
\leq&\frac{n }{3\varepsilon}\frac{|\nabla\phi|^2}{\phi}
+\frac{n\varepsilon }{3h^2}\phi|\nabla h|^2\\
=&\frac{n
}{3\varepsilon}\frac{|\nabla\phi|^2}{\phi}+\frac{n\varepsilon
}{3h^2}G,
\endaligned$$  where $\varepsilon\in(0,\frac{5}{3})$ is a positive constant, into
\eqref{Proof19} gives
\begin{equation}\label{Proof20}\aligned
\frac{(5-3\varepsilon)n }{9}\frac{G}{h^2}\leq&2\phi(a+K)-
\Delta\phi+\Big(2+\frac{n
}{3\varepsilon}\Big)\frac{|\nabla\phi|^2}{\phi}\\
\leq&2(a+K)-\Delta\phi+\Big(2+\frac{n
}{3\varepsilon}\Big)\frac{|\nabla\phi|^2}{\phi}.
\endaligned\end{equation}
In particular, choosing $\varepsilon=\frac{1}{3}$ in \eqref{Proof20}
and applying \eqref{Proof16} and \eqref{Proof17}, we obtain
\begin{equation}\label{Proof21}\aligned  \frac{4nG}{9h^2}
\leq&2(a+K)-\Delta\phi+ (n+2)
\frac{|\nabla\phi|^2}{\phi}\\
\leq&2(a+K)+\frac{C}{R^2}\left( 1+ \sqrt{K} \coth (\sqrt{K}
R)\right).
\endaligned\end{equation} It follows that for $x\in B_p(R)$,
\begin{equation}\label{Proof22}\aligned
\frac{4n}{9}G(x)\leq&\frac{4n}{9}G(x_0)\\
\leq&h^2(x_0)\left[\frac{(n+3)^2}{2n}(a+K)+\frac{C}{R^2}\left( 1+
\sqrt{K}R \coth (\sqrt{K} R)\right)\right].\endaligned\end{equation}
This shows
\begin{equation}\label{Proof23}\aligned
|\nabla
u|^2(x)\leq&M^2\left[\frac{(n+3)^2}{2n}(a+K)+\frac{C}{R^2}\left( 1+
\sqrt{K}R \coth (\sqrt{K} R)\right)\right]
\endaligned\end{equation}
and
\begin{equation}\label{Proof24}\aligned  |\nabla u|
\leq& M\sqrt{\frac{(n+3)^2}{2n}(a+K)+\frac{C}{R^2}\left( 1+
\sqrt{K}R \coth (\sqrt{K} R)\right)},
\endaligned\end{equation} where $M=\sup\limits_{x\in
B_p(2R)}u(x)$. This yields the desired inequality \eqref{1thm1} of
Theorem \ref{thm1-1}.

\vspace*{2mm}

Next, we prove the case of $a<0$.

Define $\tilde{G}=\phi |\nabla h|^2$. We apply maximum principle to
$\tilde{G}$ on $B_p(2 R)$. Assume $\tilde{G}$ achieves its maximum
at the point $\tilde{x}_0\in B_p(2R)$ and assume
$\tilde{G}(\tilde{x}_0)>0$ (otherwise the proof is trivial). Then at
the point $\tilde{x}_0$, it holds that
$$\Delta \tilde{G}\leq 0, \ \ \
\nabla (|\nabla h|^2)=-\frac{|\nabla h|^2}{\phi} \nabla \phi$$ and
\begin{equation}\label{Proof25}\aligned
0\geq& \Delta \tilde{G}\\
=&\phi\Delta(|\nabla h|^2)+\frac{\Delta\phi}{\phi}\tilde{G}
-2\frac{|\nabla\phi|^2}{\phi^2}\tilde{G}\\
\geq&2\phi\Big[\frac{37n}{36}\frac{|\nabla
h|^4}{h^2}-\frac{5n}{6}\frac{\nabla h}{h}\nabla (|\nabla
h|^2)-(a+K)|\nabla h|^2\Big]\\
&+\frac{\Delta\phi}{\phi}\tilde{G}-2\frac{|\nabla\phi|^2}{\phi^2}\tilde{G}\\
=&\frac{37n}{18}\frac{\tilde{G}^2}{\phi
h^2}+\frac{5n\tilde{G}}{3\phi}\nabla\phi\frac{\nabla
h}{h}-2(a+K)\tilde{G}\\
&+\frac{\Delta\phi}{\phi}\tilde{G}
-2\frac{|\nabla\phi|^2}{\phi^2}\tilde{G},
\endaligned\end{equation} where the second inequality  used
\eqref{Proof15}. Multiplying both sides of \eqref{Proof25} by
$\frac{\phi}{\tilde{G}}$ yields
\begin{equation}\label{Proof26}\aligned
\frac{37n}{18}\frac{\tilde{G}}{
h^2}\leq-\frac{5n}{3}\nabla\phi\frac{\nabla h}{h}+2\phi(a+K)-
\Delta\phi+2\frac{|\nabla\phi|^2}{\phi}.
\endaligned\end{equation}
Inserting the Cauchy inequality
$$\aligned-\frac{5n}{3}\nabla\phi\frac{\nabla
h}{h}\leq&\frac{5n}{3}|\nabla\phi|\frac{|\nabla
h|}{h}\\
\leq&\frac{5n
}{6\varepsilon}\frac{|\nabla\phi|^2}{\phi}+\frac{5n\varepsilon
}{6h^2}\tilde{G},
\endaligned$$  where $\varepsilon\in(0,\frac{37}{15})$ is a positive constant, into
\eqref{Proof26} gives \begin{equation}\label{Proof27}\aligned
\frac{(37-15\varepsilon)n }{18}\frac{\tilde{G}}{h^2}\leq&2\phi(a+K)-
\Delta\phi+\Big(2+\frac{5n
}{6\varepsilon}\Big)\frac{|\nabla\phi|^2}{\phi}\\
\leq&2\max\{0,\, a+K \}-\Delta\phi+\Big(2+\frac{5n
}{6\varepsilon}\Big)\frac{|\nabla\phi|^2}{\phi}.
\endaligned\end{equation}
Hence, choosing $\varepsilon=\frac{1}{15}$ in \eqref{Proof20} and
applying \eqref{Proof16} and \eqref{Proof17}, we have
\begin{equation}\label{Proof21}\aligned  \frac{2n\tilde{G}}{h^2}
\leq&2\max\{0,\, a+K \}+\frac{C}{R^2}\left( 1+ \sqrt{K} \coth
(\sqrt{K} R)\right).
\endaligned\end{equation}
Finally, restricting to $B_p(R)$ yields
\begin{equation}\label{Proof28}\aligned
|\nabla u|^2 \leq& M^2\left[\frac{(5n+6)^2}{36n}\max\{0,\, a+K
\}+\frac{C}{R^2}\left( 1+ \sqrt{K}R \coth (\sqrt{K}
R)\right)\right].
\endaligned\end{equation}
This concludes the proof of inequality \eqref{1thm2} of Theorem
\ref{thm1-1}.

\bibliographystyle{Plain}

\end{document}